# Преобразование функции Миттага-Лефлера к экспоненциальной функции и ее некоторые применения к задачам с дробной производной

**Алиев Фикрет А., Алиев Н.А., Сафарова Н.А., Гасымова К.Г.**

*Институт Прикладной Математики БГУ, Баку, Азербайджан*
*e-mail: f_aliev@yahoo.com, nihan.aliev@gmail.com, narchis2003@yahoo.com*

***Резюме:*** Впервые приводится соотношение между функцией Миттага-Лефлера к экспоненциальной. Результаты применяются к построению решения задачи Коши для обыкновенных линейных операторных дифференциальных уравнений с постоянными коэффициентами и дробными производными. На примере показывается, что когда порядок производных (дробных) приближается к целым числам, тогда результаты совпадают с классическими.

***Ключевые слова:*** экспоненциальная функция, функция Миттага-Лефлера, дробно-производная, задача Коши, линейные операторные дифференциальные уравнения.

**§1. Введение.** Несмотря на то, что в последнее время бурно развиваются разные задачи из теории управления, оптимизации, газовой динамики [1-4] и др. с дробными производными на основе функции Миттага-Лефлера до настоящего времени не установлена ее связь с экспоненциальными функциями [2]. Поэтому многие вопросы, связанные например с теорией устойчивости [5,6] остаются открытыми.

В настоящей заметке приводится соотношение между этими функциями и их результат применяется к нахождению решения задачи Коши для линейных обыкновенных операторных дифференциальных уравнений постоянными коэффициентами, с дробными производными. Далее приводится конкретный пример иллюстрирующий приближение решение к классическим решениям, когда порядок дробных производных стремится к целым числам.

**§2. Функция Миттага-Лефлера.** Функции Миттага-Лефлера представляется в виде [2,7]:

$$\mathrm{E}_\alpha(z) = \sum_{k=0}^{\infty} \frac{z^k}{\Gamma(\alpha k + 1)}, \qquad \alpha > 0 \qquad (1)$$

и

$$\mathrm{E}_{\alpha,\beta}(z) = \sum_{k=0}^{\infty} \frac{z^k}{\Gamma(\alpha k + \beta)}, \qquad \alpha > 0, \ \beta > 0, \qquad (2)$$

где $\Gamma(\cdot)$ является Гамма функцией Эйлера [2], $z$-любое действительное или комплексное переменные.

Пологая замену $z = x^\alpha$ и, учитывая формулу $\Gamma(\alpha k + 1) = (\alpha k)!$, из (1) и (2) получим:

$$E_\alpha(x^\alpha) = \sum_{k=0}^{\infty} \frac{x^{\alpha k}}{(\alpha k)!}, \qquad \alpha > 0, \qquad (3)$$

и

$$E_{\alpha,\beta}(x^\alpha) = \sum_{k=0}^{\infty} \frac{x^{\alpha k}}{(\alpha k + \beta - 1)!} = x^{1-\beta} \sum_{k=0}^{\infty} \frac{x^{\alpha k + \beta - 1}}{(\alpha k + \beta - 1)!}, \qquad (4)$$

$$\alpha > 0, \ \beta > 0.$$

Если обратить внимание на формулу из [2]

$$D^\beta \left( \frac{x^\alpha}{\alpha!} \right) = \frac{x^{\alpha - \beta}}{(\alpha - \beta)!}, \qquad \beta \le \alpha, \qquad (5)$$

и учитывая, что $\Gamma(0) = (-1)! = \infty,$ то далее имеем

$$\frac{x^{-1}}{(-1)!} = \delta(x), \qquad (6)$$

где $\delta(x)$ — дельта функция Дирака [2].

**§3. Связь функции Миттага-Лефлера с экспоненциальный функций** $e^x$. Тогда, аналог инвариантной функции Эйлера $e^x$ из аддитивного анализа [8] для производной дробного порядка получается из функции Миттага-Лефлера с помощью следующей формулы

$$h_\alpha(x) = \sum_{k=1}^{\infty} \frac{x^{-1+k\alpha}}{(-1+k\alpha)!} = \frac{x^{\alpha-1}}{(\alpha-1)!} + \frac{x^{2\alpha-1}}{(2\alpha-1)!} + \ldots . \qquad (7)$$

Действительно,

$$D^\alpha h_\alpha(x) = D^\alpha \left[ \frac{x^{\alpha-1}}{(\alpha-1)!} + \frac{x^{2\alpha-1}}{(2\alpha-1)!} + \ldots \right] = \frac{x^{-1}}{(-1)!} + \frac{x^{\alpha-1}}{(\alpha-1)!} + \frac{x^{2\alpha-1}}{(2\alpha-1)!} \ldots \qquad (7')$$

Таким образом, как следует из (6), первое слагаемое последнего является $\delta(x)$, т.е. дельта функция Дирака. А остальные являются $h_\alpha(x)$. Если $x > 0,$ то из (7') имеем:

$$D^\alpha h_\alpha(x) = h_\alpha(x). \qquad (8)$$

Поэтому построим новую функцию $h_\alpha(x,\lambda)$ являющейся аналогом функции Эйлера $e^{\lambda x}$ для дробной производной исходя из (7) в следующем виде:

$$h_\alpha(x,\lambda) = \sum_{k=1}^{\infty} \lambda^{k-1} \frac{x^{k\alpha-1}}{(k\alpha-1)!} = \frac{x^{\alpha-1}}{(\alpha-1)!} + \lambda \frac{x^{2\alpha-1}}{(2\alpha-1)!} + \lambda^2 \frac{x^{3\alpha-1}}{(3\alpha-1)!} + \ldots . \qquad (9)$$

Действительно:

$$D^{\alpha}h_{\alpha}(x,\lambda) = D^{\alpha}\sum_{k=1}^{\infty}\lambda^{k-1}\frac{x^{k\alpha-1}}{(k\alpha-1)!} = \sum_{k=1}^{\infty}\lambda^{k-1}\frac{x^{(k-1)\alpha-1}}{[(k-1)\alpha-1]!} = \frac{x^{-1}}{(-1)!} +$$

$$+ \lambda\frac{x^{\alpha-1}}{(\alpha-1)!} + .\lambda^{2}\frac{x^{2\alpha-1}}{(2\alpha-1)!} + ...,$$

Предполагая, что $x > 0$ $(\delta(x) \equiv 0)$, исходя из (6), имеем:

$$D^{\alpha}h_{\alpha}(x,\lambda) = \lambda h_{\alpha}(x,\lambda), \qquad (10)$$

Теперь, возвращаясь к анализу или алгебре [8,9], покажем, что для любого вещественного $\alpha \in (0,1]$ существуют натуральные числа $m$ и $n$ $(m,n \in N, m \leq n)$ такие, что для любого $\varepsilon > 0$

$$\left|\alpha - \frac{2m+1}{2n+1}\right| < \varepsilon, \qquad (11)$$

причем эти $m$ и $n$ не единственны.

Действительно, к каждому вещественному числу можно сколь угодно точно приближаться рациональным числом, а к любому рациональному числу можно сколь угодно точно приближаться числом вида $\frac{2m+1}{2n+1}$.

Пусть $\alpha \in (0,1), \varepsilon = 10^{-k}, k \in N$, тогда имеет место представление

$$\alpha = 0,\alpha_{1}\alpha_{2}..\alpha_{k-1}\alpha_{k}\alpha_{k+1}...,$$

где каждый $\alpha_i$ принимает одно из чисел от нуля до девяти. Тогда рассмотрим следующее рациональное число:

$$\frac{\alpha_{1}\alpha_{2}..\alpha_{k}\alpha_{k+1}\alpha_{k+2}}{10^{k+2}} \approx \frac{\alpha_{1}\alpha_{2}..\alpha_{k+2} + q}{10^{k+2} + p}$$

где $p > 2$ - число простое, а $q \in N$ такое, что $(\alpha_1,\alpha_2,...,\alpha_{k+2} + q)$ являлось бы нечетным числом.

Таким образом, существуют $m,n \in N, m < n$ такие, что

$$\frac{\alpha_{1}\alpha_{2}..\alpha_{k+2} + q}{10^{k+2} + p} = \frac{2m+1}{2n+1} \approx \alpha.$$

**Пример.** Пусть

$$\alpha = \sqrt{2} - 1 \approx 0,414 = \frac{414}{1000} \approx \frac{415}{1001} = \frac{2*207+1}{2*500+1} \approx 0,41456$$

т.е., из последнего соотношения m $= 207$, n $= 500$ получено $\varepsilon = 0,00056$.

Таким образом, задавая сколь угодно малое $\varepsilon > 0$ можно подобрать соответствующим $m$ и $n$. Этим установлено следующее утверждение.

**Лемма.** Для любого вещественного числа $\alpha \in (0,1)$ существуют такие натуральные числа $m,n \in N, m < n$, что для любого $\varepsilon > 0$ имеем

$$\left|\alpha - \frac{2m+1}{2n+1}\right| < \varepsilon.$$

Теперь рассмотрим следующую функцию

$$h_{\frac{1}{2n+1}}\left(x,\lambda^{\frac{1}{2m+1}}\right)=\sum_{k=1}^{\infty}\lambda^{\frac{k-1}{2m+1}}\frac{x^{\frac{k}{2n+1}-1}}{\left(\frac{k}{2n+1}-1\right)!}=\frac{x^{\frac{1}{2n+1}-1}}{\left(\frac{1}{2n+1}-1\right)!}+\lambda^{\frac{1}{2m+1}}\frac{x^{\frac{2}{2n+1}-1}}{\left(\frac{2}{2n+1}-1\right)!}+\ldots, \quad (12)$$

Тогда при x>0 имеем

$$D^{\frac{1}{2n+1}}h_{\frac{1}{2n+1}}\left(x,\lambda^{\frac{1}{2m+1}}\right)=\frac{x^{-1}}{(-1)!}+\lambda^{\frac{1}{2m+1}}\frac{x^{\frac{1}{2n+1}-1}}{\left(\frac{1}{2n+1}-1\right)!}+\lambda^{\frac{2}{2m+1}}\frac{x^{\frac{2}{2n+1}-1}}{\left(\frac{2}{2n+1}-1\right)!}+\ldots=$$

$$=\delta(x)+\lambda^{\frac{1}{2m+1}}\left[\frac{x^{\frac{1}{2n+1}-1}}{\left(\frac{1}{2n+1}-1\right)!}+\lambda^{\frac{1}{2m+1}}\frac{x^{\frac{2}{2n+1}-1}}{\left(\frac{2}{2n+1}-1\right)!}+\ldots\right]=\lambda^{\frac{1}{2m+1}}h_{\frac{1}{2n+1}}\left(x,\lambda^{\frac{1}{2m+1}}\right) \quad (13)$$

Далее, займемся функцией (12) и постараемся привести её к экспоненциальному виду. Для простоты рассмотрим следующую функцию:

$$h_{\frac{1}{2n+1}}(x,\rho)=\frac{x^{\frac{1}{2n+1}-1}}{\left(\frac{1}{2n+1}-1\right)!}+\rho\frac{x^{\frac{2}{2n+1}-1}}{\left(\frac{2}{2n+1}-1\right)!}+\ldots=\sum_{k=1}^{\infty}\rho^{k-1}\frac{x^{\frac{k}{2n+1}-1}}{\left(\frac{k}{2n+1}-1\right)!}.$$

Последнее выражение представим в виде суммы $(2n+1)$ ряда в следующем виде

$$h_{\frac{1}{2n+1}}(x,\rho)=\rho^{0}\frac{x^{\frac{1}{2n+1}-1}}{\left(\frac{1}{2n+1}-1\right)!}+\rho^{2n+1}\frac{x^{\frac{1}{2n+1}}}{\frac{1}{2n+1}!}+\rho^{4n+2}\frac{x^{1+\frac{1}{2n+1}}}{\left(1+\frac{1}{2n+1}\right)!}+\ldots+$$

$$+\rho^{k(2n+1)}\frac{x^{k-1+\frac{1}{2n+1}}}{\left(k-1+\frac{1}{2n+1}\right)!}+\ldots+\rho\frac{x^{\frac{2}{2n+1}-1}}{\left(\frac{2}{2n+1}-1\right)!}+\rho^{2n+2}\frac{x^{\frac{2}{2n+1}}}{\frac{2}{2n+1}!}+$$

$$+\rho^{4n+3}\frac{x^{1+\frac{2}{2n+1}}}{\left(1+\frac{2}{2n+1}\right)!}+\ldots+\rho^{k(2n+1)+1}\frac{x^{k-1+\frac{2}{2n+1}}}{\left(k-1+\frac{2}{2n+1}\right)!}+\ldots+\rho^{2}\frac{x^{\frac{3}{2n+1}-1}}{\left(\frac{3}{2n+1}-1\right)!}+$$

$$\rho^{2n+3}\frac{x^{\frac{3}{2n+1}}}{\frac{3}{2n+1}!}+\rho^{4n+4}\frac{x^{1+\frac{3}{2n+1}}}{\left(1+\frac{3}{2n+1}\right)!}+\ldots+\rho^{k(2n+1)+2}\frac{x^{k-1+\frac{3}{2n+1}}}{\left(k-1+\frac{3}{2n+1}\right)!}+\ldots+ \quad (14)$$

$$+ \rho^{2n}\frac{x^0}{0!} + \rho^{4n+1}\frac{x^1}{1!} + \rho^{6n+2}\frac{x^2}{2!} + \ldots + \rho^{k(2n+1)+2n}\frac{x^k}{k!} + \ldots =$$

$$= \sum_{k=1}^{\infty} \rho^{(k-1)(2n+1)} \frac{x^{\frac{1}{2n+1}+k-2}}{\left(\frac{1}{2n+1}+k-2\right)!} + \sum_{k=1}^{\infty} \rho^{(k-1)(2n+1)+1} \frac{x^{\frac{2}{2n+1}+k-2}}{\left(\frac{2}{2n+1}+k-2\right)!} +$$

$$+ \sum_{k=1}^{\infty} \rho^{(k-1)(2n+1)+2} \frac{x^{\frac{3}{2n+1}+k-2}}{\left(\frac{3}{2n+1}+k-2\right)!} + \ldots + \sum_{k=1}^{\infty} \rho^{(k-1)(2n+1)+2n} \frac{x^{k-1}}{(k-1)!} \equiv$$

$$\equiv J_0(x,\rho) + J_1(x,\rho) + J_2(x,\rho) + \ldots + J_{2n}(x,\rho).$$

Тогда

$$J_s(x,\rho) = \rho^s \frac{x^{\frac{s+1}{2n+1}-1}}{\left(\frac{s+1}{2n+1}-1\right)!} + \rho^{s+2n+1} \frac{x^{\frac{s+1}{2n+1}}}{\left(\frac{s+1}{2n+1}\right)!} + \ldots +$$

$$+ \rho^{s+k(2n+1)} \frac{x^{\frac{s+1}{2n+1}+k-1}}{\left(\frac{s}{2n+1}+k-1\right)!} + \ldots, \qquad s = \overline{0,2n}. \tag{15}$$

Легко видеть, что [2], интегрируя (15) с порядком $\left(1 - \frac{s+1}{2n+1}\right)$ имеем

$$I_0^{1-\frac{s+1}{2n+1}} J_s(x,\rho) = \rho^s \frac{x^0}{0!} + \rho^{s+2n+1} \frac{x^1}{1!} + \ldots + \rho^{s+k(2n+1)} \frac{x^k}{k!} + \ldots + \overline{J}_s(x,p) \quad s = \overline{0,2n}. \tag{16}$$

где функция $\overline{J}_s(x,\rho)$ подобна обычной интегральной постоянной

$$D^{1-\frac{s+1}{2n+1}} \overline{J}_s(x,\rho) = 0, \quad s = \overline{0,2n}. \tag{17}$$

Таким образом, дифференцируя (16) с порядком $\left(1 - \frac{s+1}{2n+1}\right)$ имеем

$$J_s(x,\rho) = D^{1-\frac{s+1}{2n+1}} I_0^{1-\frac{s+1}{2n+1}} J_s(x,\rho) = \rho^s D^{1-\frac{s+1}{2n+1}} e^{\rho^{2n+1}x} =$$

$$= \rho^s \frac{d}{dx} \int_0^x \frac{(x-t)^{\frac{s+1}{2n+1}-1}}{(\frac{s+1}{2n+1}-1)!} e^{\rho^{2n+1}t} dt, \; s = \overline{0,2n}. \tag{18}$$

Подставляя (18) в (14) имеем:

$$h_{\frac{1}{2n+1}}(x,\rho) = \sum_{s=0}^{2n} \rho^s \frac{d}{dx} \int_0^x \frac{(x-t)^{\frac{s+1}{2n+1}-1}}{\left(\frac{s+1}{2n+1}-1\right)!} e^{\rho^{2n+1}t} dt = \sum_{s=0}^{2n-1} \rho^s \frac{d}{dx} \int_0^x \frac{(x-t)^{\frac{s-2n}{2n+1}}}{\left(\frac{s-2n}{2n+1}\right)!} e^{\rho^{2n+1}t} dt + \rho^{2n} e^{\rho^{2n+1}x}. \tag{19}$$

Этим установлено следующее утверждение.

***Теорема1.*** Пусть $m, n \in N$, $m < n$, тогда справедливо следующая формула

$$D^{\frac{2m+1}{2n+1}} h_{\frac{1}{2n+1}}\left(x, \lambda^{\frac{1}{2m+1}}\right) = \lambda h_{\frac{1}{2n+1}}\left(x, \lambda^{\frac{1}{2m+1}}\right).$$

***Замечание.*** Легко можно видеть, что при $m = n = 0$ имеем:
$$Dh_1(x,\lambda) = \lambda h_1(x,\lambda)$$

и

$$h_1(x,\lambda) = e^{\lambda x}.$$

**§4. Решение задачи Коши для линейных обыкновенных операторных дифференциальных уравнение с дробными производными.** Иллюстрируем выше приведенные результаты при решении задачи Коши для обыкновенных операторных дифференциальных уравнений дробного порядка, т.е. имеем

$$D^{\alpha p} y + a_1 D^{\alpha(p-1)} y + \ldots + a_{p-1} D^{\alpha} y + a_p y = 0, \quad x > x_0 > 0 \qquad (20)$$

с начальными условиями

$$D^{\alpha k} y(x)\big|_{x=x_0} = \beta_k, \quad k = \overline{0, p-1}, \qquad (21)$$

где $a_k \ (k = \overline{1,p})$ и $\beta_k, \ k = \overline{0, p-1}$ действительные постоянные, $\alpha \in (0,1)$,

Рассмотрим функцию $h_\alpha(x)$ инвариантную относительно производной порядка $\alpha$, как в (7), где $\alpha$ порядок производной. Как видно из (8), функция (7) действительно является инвариантной для производной порядка $\alpha$.

Будем искать решение уравнение (20) в виде (9), тогда:
$$D^{\alpha k} h_\alpha(x,\lambda) = \lambda^k h_\alpha(x,\lambda). \qquad (22)$$

Подставляя (22) в (20) получим:
$$\left(\lambda^p + a_1 \lambda^{p-1} + \ldots + a_{p-1} \lambda + a_p\right) h_\alpha(x,\lambda) \equiv 0,$$

где для ненулевого решения уравнения (20) получаем следующее характеристическое уравнение

$$\lambda^p + a_1 \lambda^{p-1} + \ldots + a_{p-1} \lambda + a_p = 0. \qquad (23)$$

Решив характеристическое уравнение (23), находим ее корни $\lambda_1, \lambda_2, \ldots, \lambda_p$, где для простоты предположим, что они различные. Тогда общее решение уравнение (20) имеет вид:

$$y(x) = \sum_{k=1}^{p} c_k h_\alpha(x, \lambda_k), \qquad (24)$$

где $c_k \ (k = \overline{1,p})$ произвольные постоянные числа, которые определяются из начального условия (21).

Действительно, вычислив производные из (24)

$$D^{\alpha k}y(x)=\sum_{i=1}^{p}c_i\lambda_i^k h_\alpha(x,\lambda_i), \qquad (25)$$

и подставив последние в (21) получим:

$$\sum_{i=1}^{p}c_i\lambda_i^k h_\alpha(x_0,\lambda_i)=\beta_k, \quad k=\overline{0,p-1}. \qquad (26)$$

Определив произвольные постоянные $c_i$ из системы линейных алгебраический уравнений(26), в виде:

$$c_s=\frac{\Delta_s}{\Delta}, (s=\overline{1,p}), \qquad (27)$$

где

$$\Delta=\left(\prod_{i=1}^{p}h_\alpha(x_0,\lambda_i)\right)\begin{vmatrix}1 & 1 & \cdots & 1 \\ \lambda_1 & \lambda_2 & \cdots & \lambda_p \\ \cdots & \cdots & \cdots & \cdots & \cdots \\ \lambda_1^{p-1} & \lambda_2^{p-1} & \cdots & \lambda_p^{p-1}\end{vmatrix}, \qquad (28)$$

$$\Delta_s=\left(\prod_{\substack{i=1\\i\neq s}}^{p}h_\alpha(x_0,\lambda_i)\right)\begin{vmatrix}1 & 1\beta_0 1 & \ldots & 1 \\ \lambda_i & \lambda_{s-1}\beta_1\lambda_{s+1} & \ldots & \lambda_p \\ \cdots & \cdots & \cdots & \cdots & \cdots & \cdots \\ \lambda_1^{p-1} & \lambda_{s-1}^{p-1}\beta_{p-1}\lambda_{s+1}^{p-1} & \ldots & \lambda_p^{p-1}\end{vmatrix}, \quad s=\overline{1,p}. \qquad (29)$$

Переходя от функций Миттага-Лефлер (9), (14) к экспоненциальным функций (19), легко можно видеть, что выражение для решения граничный задачи (20), (21) как видно из (24)-(27), примет вид

$$y(x)=\sum_{k=1}^{p}\frac{\Delta_k}{\Delta}\left\{\sum_{s=0}^{2n-1}\lambda_k^{\frac{s}{2m+1}}\frac{d}{dx}\int_0^x\frac{(x-t)^{\frac{s-2n}{2n+1}}}{\left(\frac{s-2n}{2n+1}\right)!}e^{\lambda_k^{\frac{2n+1}{2m+1}}t}dt+\lambda_k^{\frac{2n}{2m+1}}e^{\lambda_k^{\frac{2n+1}{2m+1}}x}\right\}. \qquad (30)$$

Таким образом, доказано следующее утверждение.

**Теорема 2:** Пусть $a_i(i=\overline{1,p}), \beta_k(k=\overline{0,p-1})$ заданные вещественные постоянные, $\alpha\in(0,1)$ тогда, задача Коши (20), (21) при $x\geq x_0>0$ имеет единственное решение, которое представимо в виде (30), где $\Delta$ и $\Delta_s$ определяются из (28) и (29) соответственно.

**Замечание:** При $\alpha=1$ задача (20), (21) превращается в классическую задачу Коши [10] для уравнения $p$-го порядка, а решение (30) при $n=0, m=0$ примет вид

$$y(x)=\sum_{k=1}^{p}\frac{\Delta_k}{\Delta}e^{\lambda_k x}.$$

Как видно из (28) и (29):

$$\Delta = \prod_{i=1}^{p} e^{\lambda_i x_0} \begin{vmatrix} 1 & 1 & \cdots & 1 \\ \lambda_1 & \lambda_2 & \cdots & \lambda_p \\ \cdots & \cdots & \cdots & \cdots \\ \lambda_1^{p-1} & \lambda_2^{p-1} & \cdots & \lambda_p^{p-1} \end{vmatrix},$$

а $\Delta_s$ получается подобным образом.

Теперь покажем, что при $x \to \infty$ решение (30) стремится к нулю, когда $\lambda_i < 0, (i = \overline{1, p})$, т.е. $y(x) \to 0$ при $x \to \infty$. Тогда преобразование (30) на следующем виде

$$\frac{d}{dx}\int_0^x \frac{(x-t)^{\frac{s-2n}{2n+1}}}{\frac{s-2n}{2n+1}!} e^{\lambda k^{\frac{2n+1}{2m+1}}t} dt = -\frac{d}{dx}\int_0^x \frac{d}{dt}\left(\frac{(x-t)^{\frac{s+1}{2n+1}}}{\frac{s+1}{2n+1}!}\right) e^{\lambda k^{\frac{2n+1}{2m+1}}t} dt =$$

$$= -\frac{d}{dx}\left[\frac{(x-t)^{\frac{s+1}{2n+1}}}{\frac{s+1}{2n+1}!} e^{\lambda k^{\frac{2n+1}{2m+1}}t}\bigg|_{t=0}^{x} - \int_0^x \frac{(x-t)^{\frac{s+1}{2n+1}}}{\frac{s+1}{2n+1}!} e^{\lambda k^{\frac{2n+1}{2m+1}}t} \lambda_k^{\frac{2n+1}{2m+1}} dt\right] =$$

$$= -\frac{d}{dx}\left[-\frac{x^{\frac{s+1}{2n+1}}}{\frac{s+1}{2n+1}!} - \lambda_k^{\frac{2n+1}{2m+1}}\int_0^x \frac{(x-t)^{\frac{s+1}{2n+1}}}{\frac{s+1}{2n+1}!} e^{\lambda k^{\frac{2n+1}{2m+1}}t} dt\right] =$$

$$= \frac{x^{\frac{s-2n}{2n+1}}}{\frac{s-2n}{2n+1}!} + \lambda^{\frac{2n+1}{2m+1}}\int_0^x \frac{(x-t)^{\frac{s-2n}{2n+1}}}{\frac{s-2n}{2n+1}!} e^{\lambda k^{\frac{2n+1}{2m+1}}t} dt,$$

$$\lambda_k < 0,$$

где при $x \to \infty$ решение стремится к нулю, т.е. при $x \to \infty$ из (30) $y(x) \to 0$.

**Пример:** Пусть $p = 1$. Тогда задача Коши (20) и (21) переходит к виду

$$D^\alpha y + a_1 y = 0, \quad x \geq x_0 > 0. \qquad (31)$$

с начальным условием $y(x_0) = \beta_0$, т.е. решение (30) будет $y(x) = \frac{\beta_0}{h_\alpha(x_0, -a_1)} h_\alpha(x, -a_1)$, где с помощью экспоненциальной функции имеет вид:

$$y(x) = \beta_0 \frac{\sum_{s=0}^{2n}(-a_1)^{\frac{s}{2n+1}} \frac{d}{dx}\int_{x_0}^{x} \frac{(x-t)^{\frac{s-2n}{2n+1}}}{\left(\frac{s-2n}{2n+1}\right)!} e^{(-a_1)^{\frac{2n+1}{2m+1}}t} dt}{\sum_{s=0}^{2n}(-a_1)^{\frac{s}{2n+1}} \frac{d}{dx}\int_{x_0}^{x} \frac{(x-t)^{\frac{s-2n}{2n+1}}}{\frac{s-2n}{2n+1}!} e^{(-a_1)^{\frac{2n+1}{2m+1}}t} dt \bigg|_{x=x_0}}, \qquad (32)$$

Отметим, что при $\alpha = 1$ т.е. при $m = n = 0$ из (32) получим:

$$y(x) = \beta_0 \frac{\frac{d}{dx}\int_{x_0}^{x} e^{-a_1 t} dt}{\frac{d}{dx}\int_{x_0}^{x} e^{-a_1 t} dt\bigg|_{x=x_0}} = \beta_0 \frac{e^{-a_1 x}}{e^{-a_1 x_0}} = \beta_0 e^{-a_1(x-x_0)}.$$

Таким образом, было показано, что при $\alpha = 1$ решение (32) уравнения (31) превращается в классическую формулу [10].

Отметим, что, представление решение (30) задачи Коши для дробных производных (20),(21) при $\lambda_k < 0$ и $x \to \infty$ стремится к нулю как экспоненциальная функция. Для примера (31) из формул (32) видно, что когда $a_1 > 0$ при $t \to \infty$ решение $y(x)$ стремится к нулю.

Отметим, что можно распространит полученные результаты (30), для решения матричного случая, т.е. когда рассматривается вместо (20) системы с матричными коэффициентами [11]. А эти позволить рассмотрит задачи оптимальной стабилизации [11-15] и задачи оптимизации с неразделенными двухточечными краевыми условиями [14, 16-18].

*Заключение.* Впервые функции Миттага-Лефлера представляется с помощью экспоненциальной функцией. Это позволяет аналитически построение решения задачи Коши для линейных операторных дифференциальных уравнений дробной производной. Благодаря представлению решения с помощью экспоненциальной функции показывается асимптотической устойчивость решения задачи Коши (20), (21), когда корни характеристического уравнения (23) лежать на левой полуплоскости.

## Литература

1. Monje C.A, Chen Y.Q, Vinagre B.M, Xue D., Feliu V. Fractional-Order Systems and Controls Fundamentals and Applications, Springer, London, 2010, 414 p.

2. Samko S.G, Kilbas A.A., Marichev O.I. Fractional Integrals and Derivatives: Theory and Applications, Gordon and Breach Science publishers, Yverdon, Switzerland, 1993, 780 p.
3. Vinagre B.M, Feliu V. Optimal fractional controllers for rational order systems: a special case of the Wiener-Hopf spectral factorization method, IEEE Transactions on Automatic Control, Vol.52, No.12, 2007, pp.2385-2389.
4. Алиев Ф.А., Алиев Н.А., Сафарова Н.А., Гасымова К.Г., Раджабов М.Ф. Аналитическое конструирование регуляторов для систем с дробными производными, Proceedings of IAM, Vol.6, No.2, 2017, с.252-265.
5. Барбашин Е.А. Введение в теорию устойчивости, М., Наука, 1967, 223 с.
6. Aliev F.A., Larin V.B. Optimization of Linear Control Systems: Analytical Methods and Computational Algorithms, Gordon Breach, Amsterdam, 1998, 272 p.
7. Mittag-Leffler G. Sur la representation analytique d`une branche uniformed`une function monogene, Acta Mathematica, Vol.29, 1905, pp.101- 181.
8. Фихтенгольц Г.М. Курс дифференциального и интегрального исчисления, М.: Наука, 1966, 616 с.
9. Курош, А.Г. Курс высшей алгебры, М.: Наука, 1965, 431 с.
10. Коддингтон Э.А., Левинсон Н., Теория обыкновенных дифференциальных уравнений, Издательство иностранной литературы, Москва, 1958, 475 с.
11. Андреев Ю.И. Управление конечномерными линейными объектами М., Наука ,1976, 424 с.
12. Алиев Ф.А., Ларин В.Б. Синтез оптимальных импульсивных регуляторов при идеальном измерении координат объекта. В книге Навигационные Гироскопические Системы. К.,1973,с.5-28.
13. Алиев Ф.А., Ларин В.Б. Синтез дискретных инвариантных во времени систем стабилизации. В ки.: Дискретные Системы Управления, К., 1974, с.3-22.
14. Алиев Ф.А. Методы решения прикладных задач оптимизации динамических систем, Элм, Баку ,1989, 320 с.
15. Алиев Ф.А., Бордюг Б.А., Ларин В.Б. $H_2$-оптимизация и метод пространства состояний в задаче синтеза оптимальных регуляторов , Элм, Баку, 1991, 373 с.
16. Алиев Ф.А. Выбор программы движения шагающего аппарата с почти невесомыми ногами, В кн.: Системы Навигации и Управление, Киев: Ин-т Математики АН, УССР, 1983, с.85-97.
17. Алиев Ф.А. Задачи оптимизации с двухточечными краевыми условиями. Изв. Ан СССР, техн. Кибернетика, 1985, №6,с.138-146.
18. Алиев Ф.А. Задача оптимального управления линейной системой с неразделенными граничными условиями, Дифференциальные уравнения 1986, №2, с.345-347.